\documentclass[11pt]{amsart}
\usepackage{amsthm}

\theoremstyle{plain}
\newtheorem{thm}{Theorem}[section]
\newtheorem{theorem}[thm]{Theorem}

\theoremstyle{definition}

\numberwithin{equation}{section}

 \title[Coordinate-wise Armijo's condition: General case]{Coordinate-wise Armijo's condition: General case}
 \author{Tuyen Trung Truong}
   \address{Department of Mathematics, University of Oslo, Blindern 0851 Oslo, Norway}
  \email{tuyentt@math.uio.no}

    \date{\today}
   
   \subjclass[2010]{}

\begin{document}
\maketitle

\begin{abstract}
Let $z=(x,y)$ be coordinates for the product space $\mathbb{R}^{m_1}\times \mathbb{R}^{m_2}$. Let $f:\mathbb{R}^{m_1}\times \mathbb{R}^{m_2}\rightarrow \mathbb{R}$ be a $C^1$ function, and $\nabla f=(\partial _xf,\partial _yf)$ its gradient. Fix $0<\alpha <1$.  For a point $(x,y) \in \mathbb{R}^{m_1}\times \mathbb{R}^{m_2}$, a number $\delta >0$ satisfies Armijo's condition at $(x,y)$ if the following inequality holds:
\begin{eqnarray*}
f(x-\delta \partial _xf,y-\delta \partial _yf)-f(x,y)\leq -\alpha \delta (||\partial _xf||^2+||\partial _yf||^2).  
\end{eqnarray*}
     
In one previous paper, we proposed the following {\bf coordinate-wise} Armijo's condition. Fix again $0<\alpha <1$. A pair of positive numbers $\delta _1,\delta _2>0$ satisfies the coordinate-wise variant of Armijo's condition at $(x,y)$ if the following inequality holds: 
\begin{eqnarray*}
[f(x-\delta _1\partial _xf(x,y), y-\delta _2\partial _y f(x,y))]-[f(x,y)]\leq -\alpha (\delta _1||\partial _xf(x,y)||^2+\delta _2||\partial _yf(x,y)||^2). 
\end{eqnarray*} 
Previously we applied this condition for functions of the form $f(x,y)=f(x)+g(y)$, and proved various convergent results for them. For a general function, it is crucial - for being able to do real computations - to have a systematic algorithm for obtaining $\delta _1$ and $\delta _2$ satisfying the coordinate-wise version of Armijo's condition, much like Backtracking for the usual Armijo's condition. In this paper we propose such an algorithm, and prove according convergent results.   

We then analyse and present experimental results for some functions such as $f(x,y)=a|x|+y$ (given by Asl and Overton in connection to Wolfe's method), $f(x,y)=x^3 sin (1/x) + y^3 sin(1/y)$ and Rosenbrock's function. 
\end{abstract}


\subsection{Coordinate-wise Armijo's condition} We introduce in this subsection the coordinate-wise Armijo's condition, following \cite{truong2}. For current status of results in Backtracking Gradient Descent (GD) and modifications, readers can see  \cite{truong2} as well as references therein. For feasibility and good performance of implementing Backtracking GD in Deep Neural Networks (DNN), together with source codes, see \cite{truong-nguyen} (and a more recent work in \cite{vaswani-etal}). For a review of some other popular gradient descent methods, the readers can see \cite{ruder}. 

Gradient Descent (GD) methods, invented by Cauchy in 1847 \cite{cauchy}, aim to find minima of a $C^1$ function $f:\mathbb{R}^m\rightarrow \mathbb{R}$ by the following iterative procedure
\begin{eqnarray*}
z_{n+1}=z_n-\delta (z_n)\nabla f(z_n),
\end{eqnarray*}
where $\delta (z_n)>0$, the learning rate, must be appropriately chosen. The most basic and known algorithm, Standard GD, uses $\delta (z_n)=\delta _0$ a constant. Armijo's condition \cite{armijo} is a well known criterion to find learning rates, which requires that
\begin{eqnarray*}
f(z-\delta (z)\nabla f(z))-f(z)\leq -\alpha \delta (z) ||\nabla f(z)||^2,
\end{eqnarray*}
where $0<\alpha <1$ is a given constant. 

One can implement Armijo's condition for practical computations by the following algorithm. 

{\bf Backtracking GD.} Armijo's condition gives rise to Backtracking GD, which is the following procedure. Given $0<\alpha ,\beta <1$ and $\delta _0>1$. For each $z\in \mathbb{R}^m$ we define $\delta (z)$ to be the largest number among $\{\beta ^n\delta _0:~n=0,1,2,\ldots \}$ which satisfies Armijo's condition. For each initial point $z_0$, we then define inductively the sequence $z_{n+1}=z_n-\delta (z_n)\nabla f(z_n)$. 

The purpose of the current paper is to improve Armijo's condition to adapt better to cost functions $f$, such as $a|x|+y$ or $x^3 \sin (1/x) +y^3\sin (1/y)$, where at a point the partial derivatives in different directions can be very much different in sizes. From now on, we consider the following case: $\mathbb{R}^m=\mathbb{R}^{m_1}\times \mathbb{R}^{m_2}$ is a product space with coordinates $z=(x,y)$. Extensions to more general cases such as $\mathbb{R}^{m_1}\times \mathbb{R}^{m_2}\times \mathbb{R}^{m_3}$ are straight forward. We have the following definition.

{\bf Coordinate-wise Armijo's condition.} Fix $0<\alpha <1$. A pair $\delta _1,\delta _2>0$ satisfies Coordinate-wise Armijo's condition at $z=(x,y)$ if 
\begin{eqnarray*}
f(x-\delta _1\partial _xf(x,y), y-\delta _2\partial _y f(x,y))-f(x,y)\leq -\alpha (\delta _1||\partial _xf(x,y)||^2+\delta _2||\partial _yf(x,y)||^2). 
\end{eqnarray*} 

Unlike the case of the usual Armijo's condition, it is not easy for one to choose $\delta _1$ and $\delta _2$ in such a way which makes actual computations amenable. In \cite{truong2}, we treated the special case where $f(x,y)=f_1(x)+f_2(y)$, in which case one can use the usual Backtracking algorithm separately for $f_1$ and $f_2$. In case there are cross terms between $x$ and $y$ in the function $f$, such a simple approach of choosing $\delta _1$ and $\delta _2$ are unachievable. In this paper we extend this  to  the general setting. 

\subsection{Coordinate-wise Backtracking GD}

We assume that $f$ is  a $C^1$ function on $\mathbb{R}^m\backslash A$, where $A$ is a closed subset. We fix $0<\alpha ,\beta <1$ and $\delta _0>0$. Let $z=(x,y)\in \mathbb{R}^m\backslash A$ and $r(z)=dist (z,A)>0$.  We will now construct a pair $0<\delta _x(z), \delta _y(z)<r(z)/||\nabla f(z)||$ and belonging to the discrete set $\{\beta ^n\delta _0:~n=0,1,2\ldots \}$ - in a deterministic manner - so that the coordinate-wise Armijo's condition is satisfied. 

First, we choose $\delta (z)$ to be the learning rate chosen from the usual Armijo's condition, that is $0<\delta (z)<r(z)/||\nabla f(z)||$ is the largest number $\delta $ among $\{\beta ^n\delta _0:~n=0,1,2\ldots  \}$ so that 
\begin{eqnarray*}
f(z-\delta (z)\nabla f(z))-f(z)\leq -\delta (z)||\nabla f(z)||^2. 
\end{eqnarray*}  
We recall that $||\nabla f(z)||^2=||\partial _xf(z)||^2+||\partial _yf(z)||^2$. Now we proceed in two steps: the first step is to construct $\delta _x(z)$, and then the second step is to construct $\delta _y(z)$. 

{\bf Step 1:} We choose $\delta (z)\leq \delta _x(z)<r(z)/||\nabla f(z)||$ to be the largest number $\delta _1$ among $\{\beta ^n\delta _0:~n=0,1,2\ldots  \}$ so that 
\begin{eqnarray*}
f(x-\delta _1\partial _xf(z),y-\delta (z)\partial _yf(z))-f(z)\leq -\alpha (\delta _1||\partial _xf(z)||^2+\delta (z)||\partial _yf(z)||^2).
\end{eqnarray*}
Note that there is at least one such number ($\delta _1=\delta (z)$). 

{\bf Step 2:} We choose $\delta (z)\leq \delta _y(z)< r(z)/||\nabla f(z)||$ to be the largest number $\delta _2$ among $\{\beta ^n\delta _0:~n=0,1,2\ldots  \}$ so that 
\begin{eqnarray*}
f(x-\delta _x(z)\partial _xf(z),y-\delta _2\partial _yf(z))-f(z)\leq -\alpha (\delta _x(z)||\partial _xf(z)||^2+\delta _2||\partial _yf(z)||^2).
\end{eqnarray*}
Again, note that there is at least one such number ($\delta _2=\delta (z)$). 

{\bf Remark.} If we change the order, that is to construct $\delta _y(z)$ first and then construct $\delta _x(z)$ second, then the values we obtain may be different. The following is a good heuristic to proceed. Assume that $\partial _xf$ and $\partial _yf$ are locally Lipschitz continuous near $z$, with corresponding Lipschitz constants $L_x(z)$ and $L_y(z)$ so that $L_x(z)\geq L_y(z)$, then we construct $\delta _x(z)$ first and then $\delta _y(z)$. In the opposite case of $L_x(z)\leq L_y(z)$, we switch to constructing $\delta _y(z)$ first and then $\delta _x(z)$.   

\subsubsection{Convergence analysis} First, in the case where $f$ is $C^1$ on the whole $\mathbb{R}^m$, then from $\delta _0\geq \delta _x(z),\delta _y(z)\geq \delta (z)$, where $\delta (z)$ is the learning rate from the usual Backtracking GD algorithm, we obtain the following result, whose proof is the same as in \cite{truong-nguyen}.

\begin{theorem}
Let $f$ be $C^1$ on the whole $\mathbb{R}^m$ and $\{z_n\}$ be a sequence constructed from the Coordinate-wise Backtracking GD. Then:

1) Every cluster point of $\{z_n\}$ is a critical point of $f$. 

2) Either $\lim _{n\rightarrow\infty}f(z_n)=-\infty$ or $\lim _{n\rightarrow\infty}||z_{n+1}-z_n||=0$. 

3) Let $B$ be a compact component of the set of critical points of $f$, and let $C$ be the set of cluster points of $\{z_n\}$. If $B\cap C\not= \emptyset$, then $C\subset B$ and C is connected. 

\label{Theorem1}\end{theorem}

For the special case where $f(x,y)=f_1(x)+f_2(y)$, then the numbers $\delta _x(z)$ and $\delta _y(z)$ are the same as those obtained by using Backtracking GD separately on $f_1$ and $f_2$, and hence we recover the constructions and results in \cite{truong2}.

\subsection{Some examples}
\subsubsection{Example by Asl and Overton} In \cite{asl-overton}, the following function was considered $f(x,y)=a|x|+y$, where $a>0$ is big enough. They showed that if we apply Wolfe's method with an initial point $z_0=(x_0,y_0)$ with $x_0\not=0$, then the sequence $z_n$ will converge. They argued that this showed that Wolfe's method does not work well with this example, since this function is not bounded from below, and hence if a method works well, it should not converge but should diverge. 

Here, we present $5$ small observations: 

- First, the function $g(x)=|x|$ is problematic for Standard GD and Newton's method, as well as some other familiar methods. Indeed, if we apply Standard GD to $g(x)$, for a random initial point $x_0$, then after a finite number of steps we will have that $\delta _0>2|x_n|$. After that, then the sequence $\{x_n\}$ will become periodic, going back and forth between two points. Since $f"(x)=0$ identically on $\mathbb{R}\backslash \{0\}$, it follows that Newton's method is not applicable. If we try to rescue by declaring that $f"(x)=\epsilon \not=0$ whenever $f"(x)=0$, then the sequence $\{x_n\}$ will diverge to infinity. 

- Second, Backtracking GD works well for the function $g(x)=|x|$.  Indeed, since the function $g(x)$ is coercive, it follows that the sequence $\{x_n\}$ constructed by Backtracking GD, being actually descent, is bounded. Let $x_{\infty}$ be a cluster point of $\{x_n\}$. If $x_{\infty}\not=0$, then since the function $g$ is $C^1$ in $\mathbb{R}\backslash \{0\}$, it follows from theoretical properties of Backtracking GD that $x_{\infty}$ must be a critical point of $f$. However, $f$ has no critical point inside $\mathbb{R}\backslash \{0\}$. Therefore, $x_{\infty}=0$, which means that $\{x_n\}$ converges to $0$. 

- Third, the function $f(x,y)=a|x|+y$ is problematic for Backtracking GD, in the sense given by Asl and Overton. Here the argument is not rigorous, but rather based on experiments. Indeed, in some experiments, we found that if we apply Backtracking GD to the function $f(x,y)=a|x|+y$, for an initial point $z_0=(x_0,y_0)$, then the sequence $\{z_n\}$ seems to converge. Here we offer one explanation. The number $\delta (z_n)$ so that  Armijo's condition $f(z_n-\delta (z_n))-f(z_n)\leq -\alpha ||\nabla f(z_n)||^2$ tends to depend more on the variable $x_n$ than the variable $y_n$, when $n$ large enough. More specifically, we have $\delta (z_n)\sim \delta (x_n)$, where $\delta (x_n)$ is the learning rate obtained from Backtracking GD for the function $g(x)=|x|$. Now, it seems that there is a tendency for the sequence $\{x_n\}$ to keep on {\bf the same side} (that is, either always positive or always negative, at least when $n$ is large enough). That requires that $\delta (x_n)\leq |x_n|$. However, in the Backtracking GD procedure, we will choose $\delta (x_n)$ in the discrete set $\{\delta _0, \beta \delta _0, \beta ^2\delta _0, \ldots \}$. Hence, the sequence $\{\delta (z_n)\}$ is bounded by a geometric series, and hence $\{z_n\}$ converges. 

- Fourth, the function $f(x,y)=a|x|+y$  is {\bf not problematic} for Coordinate-wise Backtracking GD, in the sense of Asl and Overton. In fact, in this case while we still have $\delta _x(z_n)\sim \delta (x_n)$, we do have $\delta _y(z_n)\sim \delta _0$. Therefore, the sequence $\{z_n\}$ diverges to infinity, more precisely to $(0,-\infty )$. 

- Fifth, from another viewpoint (optimisation on manifolds), the function $f(x,y)=a|x|+y$ is {\bf not problematic} for Backtracking GD. Indeed, let $X=\mathbb{R}^2\backslash \{x=0\}$, which is exactly the set where $f$ is differentiable. Then $X$ is a manifold, and $f:X\rightarrow \mathbb{R}$ is $C^1$ (indeed, $C^2$). Then the analog of Theorem \ref{Theorem1} in this setting is as follows: Let $\{x_n\}$ be a sequence constructed from Backtracking GD. Then (because $f$ has no critical point inside $X$) the sequence $\{x_n\}$ diverges to the boundary of $X$ (as a subset of the real projective plane $\mathbb{P}^2$).

{\bf Remark:} On the other hand, there is no problem with applying all the above numerical methods with the similar function where we replace $|x|$ by $ReLU(x)=\max\{x,0\}$. This may be one reason why Deep Neural Networks work well as observed in practice. 

\subsubsection{A singular example} Now we look at the case $f(x,y)=x^3\sin (1/x) +y^3\sin (1/y)$, which was also discussed in \cite{truong2}. 

First, we look at the one variable function $g(x)=x^3\sin (1/x)$. This function is $C^1$, but its gradient is not locally Lipschitz continuous at $x=0$. The point $x=0$ is a singular critical point of $f$: in every neighbourhood of $0$, there are both local minima and local maxima of $f$. 

If we apply Standard GD to $g(x)$, it is evident that the constructed sequence $\{x_n\}$ converges to $0$. If we start from an initial point near $0.55134554$, then after 381 iterates the point $x_n$ is $2e-09$.   

On the other hand, Backtracking GD works just fine for this function: with the same initial point, after 20 iterates the point $x_n$ is $0.24520926$. It is evident that the sequence $\{x_n\}$ converges to a local minimum near $0.24520926$.   

Second, we observe that Backtracking GD and Coordinate-wise Backtracking GD have similar performances for the function $f(x,y)$. 

{\bf Remark:} While currently we can prove avoidance of saddle points only for the modification Backtracking GD - New \cite{truong}, it appears from the experiments that Backtracking GD itself can also avoid saddle points. 

\subsubsection{Rosenbrock's function} Here we explore the performance of Coordinate-wise Backtracking GD on a truly variable-crossing landmark function. The function is $f(x,y)=(x-1)^2+100(y-x^2)^2$. This function has only one critical point at $(1,1)$, which is also a global minimum. The partial derivatives in the $x$-direction and in the $y$-direction behave very differently in this example. In this case, it happens that the order we do in each step of Coordinate-wise Backtracking GD does  matter. For a specific experiment, with initial point $(0.55134554, 0.75134554)$ we find that: 

- Standard GD encounters overflow error. 

- Backtracking GD needs 2433 iterations to reach very close to the point $(1,1)$. (The result reported in the computer is $(1.,1.)$.)

- Coordinate-wise Backtracking GD, when we construct in each step $\delta _x(z_n)$ first and then $\delta _y(z_n)$ second, needs  13342 iterations to reach very close to the point $(1,1)$. 

- On the other hand, Coordinate-wise Backtracking GD, when we construct in each step $\delta _y(z_n)$ first and $\delta _x(z_n)$ second, needs  only 4553 iterations to reach very close to the point $(1,1)$. 

\subsection{Conclusions} 

In this paper, we introduced Coordinate-wise Backtracking GD for a general $C^1$ function, generalising our previous work \cite{truong2} (where only the special case $f(x,y)=f_1(x)+f_2(y)$ was treated). The update rule is coordinate-wise based: $(x_{n+1},y_{n+1})$ $=$ $(x_n,y_n)-(\delta _x(z_n)\partial _xf(z_n),\delta _y(z_n)\partial _yf(z_n))$, where $z_n=(x_n,y_n)$ so that Armijo's condition is satisfied. Here, there is an asymmetry:  At each step, we need to make a decision on whether choose $\delta _x(z_n)$ first and then $\delta _y(z_n)$ second, or vice versa. This is to adapt better to the cases where $\partial _xf$ and $\partial _yf$ can be very different.  We prove  convergence results similar to those in \cite{truong-nguyen, truong, truong2}. 

We demonstrated that  this new algorithm allows to resolve the problem alleged in \cite{asl-overton} for the function $a|x|+y$, while Backtracking GD also is problematic for the function.  We argue that if one view optimisation of this function $a|x|+y$ as on the manifold $X=\mathbb{R}^2\backslash \{x=0\}$, then Backtracking GD itself works without problem, when one interprets the convergence results as follows: the sequence either converges to a critical point of $f$ inside $X$, or diverges to the boundary $\partial X$.  

We tested with  experiments and found that for the function $x^3\sin (1/x) +y^3\sin (1/y)$, Backtracking GD and Coordinate-wise Backtracking GD behave similarly. On the other hand, for the Rosenbrock's function, Backtracking GD and Coordinate-wise Backtracking GD behave very differently. Also, for the same function, the decision to choose either $\delta _x(z_n)$ first or $\delta _y(z_n)$  first leads to different performances. We reiterate the impression that Backtracking GD (and not just the modification Backtracking GD - New in \cite{truong}) can avoid saddle points. 

In view of performance of numerical optimization methods, there is a clear difference between the functions $|x|$ and ReLU(x), even though they have similar shapes. The numerical methods behave better for ReLU(x), and this may be a reason for why ReLU(x) is good for Deep Neural Networks as observed in reality.

\subsection{Acknowledgments} Some ideas of the work were initiated in our visit to Torus Actions SAS (Toulouse, France). We thank them for inspiring discussions and hospitality, and thank Trond Mohn Foundation for a travel fund. This work is supported by Young Research Talents grant number 300814 from Research Council of Norway.

\end{document}